\documentclass[12pt, letterpaper]{amsart}

\usepackage[english]{babel}
\usepackage{amsmath,enumerate, amsfonts, amssymb,amsthm, xypic}
\usepackage[dvips]{graphics} 

\newtheorem{thm}{Theorem}[section]

\theoremstyle{definition}

\DeclareMathOperator{\R}{\mathbf R}

\DeclareMathOperator{\Z}{\mathbf Z}

\DeclareMathOperator{\id}{id}
\DeclareMathOperator{\codim}{codim}

\title{Triangulations of \\
non-proper semialgebraic Thom maps}
\author{Masahiro Shiota}
\address{Graduate School of Mathematics, Nagoya University, Chikusa, Nagoya, 
464-8602, Japan}
\email
{shiota@math.nagoya-u.ac.jp}


\begin{document}

\begin{abstract} In \cite{S3} I solved the Thom's conjecture that a proper Thom map is triangulable. In this paper I drop the properness condition in the semialgebraic case and, moreover, in the definable case in an o-minimal structure. 
\end{abstract}

\maketitle
\section{Introduction}
Let $r$ be always a positive integer or $\infty$, $X$ and $Y$ subsets of $\R^m$ and $\R^n$, respectively, and $f:X\to Y$ a $C^r$ map (i.e., $f$ is extended to a $C^r$ map from an open neighborhood of $X$ in $\R^m$ to one of $Y$ in $\R^n$). 
A {\it $C^r$ stratification} of $f$ is a pair of $C^r$ stratifications $\{X_i\}$ of $X$ and $\{Y_j\}$ of $Y$ such that for each $i$, the image $f(X_i)$ is included in some $Y_j$ and the restriction map $f|_{X_i}:X_i\to Y_j$ is a $C^r$ submersion. 
We call also $f:\{X_i\}\to\{Y_j\}$ a {\it $C^r$ stratification} of $f:X\to Y$. 
We call $f:X\to Y$ a {\it Thom $C^r$ map} if there exists a Whitney $C^r$ stratification $f:\{X_i\}\to\{Y_j\}$ such that the following condition is satisfied. 
Let $X_i$ and $X_{i'}$ be strata with $X_{i'}\cap(\overline X_i-X_i)=\emptyset$. 
If $\{a_k\}$ is a sequence of points in $X_i$ converging to a point $b$ of $X_{i'}$ and if the sequence of the tangent spaces $\{T_{a_k}(f|_{X_i})^{-1}(f(a_k))\}$ converges to a space $T\subset\R^m$ in the Grassmannien space $G_{m,m'},\ m'=\dim (f|_{X_i})^{-1}(f(a_k))$, then $T_b(f|_{X_{i'}})^{-1}(f(b))\subset T$. 
We call then $f:\{X_i\}\to\{Y_j\}$ a {\it Thom $C^r$ stratification} of $f:X\to Y$. 
In \cite{S3} I solved the following Thom's conjecture. 

\begin{thm}
Assume $X$ and $Y$ are closed in $\R^m$ and $\R^n$, respectively, and $f:X\to Y$ is a proper Thom $C^\infty$ map. 
Then there exist homeomorphisms $\tau$ and $\pi$ from $X$ and $Y$ to polyhedra $P$ and $Q$, respectively, such that $\pi\circ f\circ\tau^{-1}:P\to Q$ is piecewise-linear. 
\end{thm}

Here a natural question arises. 
Whether can we drop the properness condition? 
Indeed, the condition is too strong for some applications. 
For example, the natural map from a $G$-manifold $M$ to its orbit space is a Thom map but not necessarily proper provided the action $G\times M\ni(g,x)\to(g x,x)\in M^2$ is proper (see \cite{MS}). 
In the present paper we give a positive answer in the semialgebraic or definable case. 
A $C^r$ stratification $f:\{X_i\}\to\{Y_j\}$ of $f:X\to Y$ is called {\it semialgebraic} ({\it definable}) if $X,\,Y,\,f,\,X_i$ and $Y_j$ are all semialgebraic (definable, respectively,) and $\{X_i\}$ and $\{Y_j\}$ are finite stratifications. 

\begin{thm}\label{thmA}
Assume $X$ and $Y$ are closed and semialgebraic (definable in an o-minimal structure) in $\R^m$ and $\R^n$, respectively, and $f:X\to Y$ is a semialgebraic (definable, respectively,) Thom $C^1$ map. 
Then there exist finite simplicial complexes $K$ and $L$ and semialgebraic (definable, respectively,) $C^0$ imbeddings $\tau:X\to|K|$ and $\pi:Y\to|L|$ such that $\tau(X)$ and $\pi(Y)$ are unions of some open simplexes of $K$ and $L$, respectively, and $\pi\circ f\circ\tau^{-1}:\tau(X)\to\pi(Y)$ is extended to a simplicial map from $K$ to $L$, where $|K|$ denotes the underlying polyhedron to $K$. 
\end{thm}

The theorem does not necessarily hold without the condition that $X$ is closed in $\R^m$. 
A counter-example is given by $X=\R^2-\{(x,y)\in\R^2:x=0,\ y\not=0\},\ Y=\R^2$ and $f(x,y)=(x,x y)$. 
Such $f$ is not triangulable in the weak sense that there exist $C^0$ imbeddings $\tau$ of $X$ and $\pi$ of $Y$ into some Euclidean space $\R^n$ such that $\overline{\tau(X)}$ is a polyhedron and $\pi\circ f\circ\tau^{-1}:\tau(Y)\to\pi(X)$ is extended to a piecewise-linear map $\theta:\overline{\tau(X)}\to\R^n$ for the following reason. 
Assume there exist $\tau$ and $\pi$ as required. 
Then $\overline{\tau(X)}$ is of dimension two and $\theta^{-1}(y)$ is of dimension 0 for each $y\in\overline{\pi(Y)}$ because $\theta$ is piecewise-linear and $\theta|_{\tau(X)}$ is injective. 
Hence a small compact neighborhood $U$ of $\tau(0)$ in $\overline{\tau(X)}$ does not intersect with $\theta^{-1}(\pi(0))$ except at $\tau(0)$. 
Choose a point $(x_1,x_2)$ in $X$ with $x_2\not=0$ so close to 0 that the half-open segment $L$ with ends $(0,x_2)$ and $(x_1,x_2)$ in $X$ is included in $\tau^{-1}(U)$. 
Then $\overline{f(L)}-f(L)=\{0\}$ and $\overline{\pi\circ f(L)}-\pi\circ f(L)=\{\pi(0)\}$. 
Hence $(\overline{\tau(L)}-\tau(L))\cap U=\{\tau(0)\}$ or $(\overline{\tau(L)}-\tau(L))\cap U=\emptyset$ since $\theta^{-1}(\pi(0))=\{\tau(0)\}$ in $U$. 
The former case contradicts the definition of $L$ and the fact that $\tau$ is a $C^0$ imbedding, and the latter does the fact that $U$ is compact. \par
  An open problem is whether a Thom $C^1$ map $f:X\to Y$ is triangulable in this weak sense under the condition that $X$ is closed in $\R^n$ or, equivalently, $X$ is locally compact. 

\section{Tube systems}
If $r$ is larger than one, $C^r$ {\it tube} at a $C^r$ submanifold $M$ of $\R^n$ is a triple $T=(|T|,\pi,\rho)$, where $|T|$ is an open neighborhood of $M$ in $\R^n$, $\pi:|T|\to M$ is a submersive $C^r$ retraction and $\rho$ is a non-negative $C^r$ function on $|T|$ such that $\rho^{-1}(0)=M$ and each point $x$ on $M$ is a unique and non-degenerate critical point of $\rho|_{\pi^{-1}(x)}$. 
We will need to consider a $C^1$ tube. 
Assume $M$ is a $C^1$ submanifold of $\R^n$. 
Let $|T|$ be an open neighborhood of $M$ in $\R^n$, $\pi:|T|\to M$ a $C^1$ map and $\rho$ a $C^1$ function on $|T|$. 
We call $T=(|T|,\pi,\rho)$ a $C^1$ {\it tube} at $M$ if there exists a $C^1$ imbedding $\tau$ of $|T|$ into $\R^n$ such that $\tau(M)$ is a $C^2$ submanifold of $\R^n$ and $\tau_*T=(\tau(|T|),\tau\circ\pi\circ\tau^{-1},\rho\circ\tau^{-1})$ is a $C^2$ tube at $\tau(M)$. 
(See pages 33--40 in \cite{S2}, which says the arguments on tube systems in \cite{G} work in the $C^1$ category.) 
A $C^r$ {\it tube system} $\{T_j\}$ for a $C^r$ stratification $\{Y_j\}$ of a set $Y\subset\R^n$ consists of one tube $T_j$ at each $Y_j$. 
We define a $C^r$ {\it weak tube system} $\{T_j=(|T_j|,\pi_j,\rho_j)\}$ for the same $\{Y_j\}$ weakening the conditions on $\rho_j$ as follows. 
Each $\rho_j$ is a non-negative $C^0$ function on $|T_j|$ with zero set $Y_j$, of class $C^r$ on $|T_j|-Y_j$ and regular on $Y_{j'}\cap\pi^{-1}_j(y)-Y_j$ for each $y\in Y_j$ and $Y_{j'}$. 
Note a $C^r$ tube system is a $C^r$ weak tube system if $\{Y_j\}$ is a Whitney stratification by Lemma I.1.1, \cite{S2}. 
In the following arguments we shrink $|T_j|$ many times without mention. \par
  We call a $C^r$ (weak) tube system $\{T_j\}$ for $\{Y_j\}$ {\it controlled} if for each pair $j$ and $j'$ with $(\overline{Y_{j'}}-Y_{j'})\cap Y_j\not=\emptyset$, 
$$
\pi_j\circ\pi_{j'}=\pi_j\quad\text{and}\quad \rho_j\circ\pi_{j'}=\rho_j\quad\text{on}\ |T_j|\cap|T_{j'}|. 
$$
Remember there exists a controlled $C^r$ tube system for a Whitney stratification (see \cite{G} and \cite{S2}), note if $\{T_j\}$ is such a $C^r$ tube system then the map $(\pi_j,\rho_j)|_{Y_{j'}\cap|T_j|}$ is a $C^r$ submersion into $Y_j\times\R$ because 
$$
(\pi_j,\rho_j)|_{Y_{j'}\cap|T_j|}\circ\pi_{j'}=(\pi_j,\rho_j)\quad\text{on}\ |T_j|\cap|T_{j'}|, 
$$
and if we assume only $\pi_j\circ\pi_{j'}=\pi_j$ on $|T_j|\cap|T_{j'}|$ then $\pi_j|_{Y_{j'}\cap|T_j|}$ is a $C^r$ submersion into $Y_j$. 
In the case of a $C^r$ weak tube system $(\pi_j,\rho_j)|_{Y_{j'}\cap|T_j|-Y_j}$ is a $C^1$ submersion into $Y_j\times\R$. 
Let $f:\{X_i\}\to\{Y_j\}$ be a $C^r$ stratification of a $C^r$ map $f:X\to Y$ between subsets of $\R^m$ and $\R^n$, respectively, $\{T_j^Y=(|Y_j^Y|,\pi_j^Y,\rho_j^Y)\}$ a controlled $C^r$ (weak) tube system for $\{Y_j\}$ and $\{T_i^X=(|T_i^X|,\pi_i^X,\rho_i^X)\}$ a $C^r$ (weak) tube system for $\{X_i\}$. 
We call $\{T_i^X\}$ {\it controlled over} $\{T_j^Y\}$ if the following four conditions are satisfied. 
Let $f$ be extended to a $C^r$ map $\tilde f:\cup_i|T_i^X|\to\R^n$. \par\noindent
(1) For each $(i,j)$ with $f(X_i)\subset Y_j$, 
$$
f\circ\pi_i^X=\pi_j^Y\circ\tilde f\quad\text{on}\ |T_j^X|\cap\tilde f^{-1}(|T_j^Y|). 
$$\par\noindent
(2) For each $j$, $\{T_i^X:f(X_i)\subset Y_j\}$ is a controlled $C^r$ (weak) tube system for $\{X_i:f(X_i)\subset Y_j\}$. \par\noindent
(3) For each pair $i$ and $i'$ with $(\overline{X_{i'}}-X_{i'})\cap X_i\not=\emptyset$,
$$
\pi_i^X\circ\pi_{i'}^X=\pi_i^X\quad\text{on}\ |T_i^X|\cap|T_{i'}^X|. 
$$\par\noindent
(4) For each $(i,j)$ with $f(X_i)\subset Y_j$ and $(i',j')$ with $(\overline{X_{i'}}-X_{i'})\cap X_i\not=\emptyset$ and $f(X_{i'})\subset Y_{j'}$, $(\pi_i^X,f)|_{X_{i'}\cap|T_i^X|}$ is a $C^r$ submersion into the fiber product $X_i\times_{(f,\pi_j^Y)}(Y_{j'}\cap|T_j^Y|)$---the $C^r$ manifold $\{(x,y)\in X_i\times(Y_{j'}\cap|T_j^Y|):f(x)=\pi_j^Y(y)\}$. \par
  Note (4) is equivalent to the next condition. \par\noindent
(4)$'$ For $(i,j),\ (i',j')$ as in (4) and for each $x\in X_{i'}\cap|T_i^X|$, the germ of $\pi_i^X|_{X_{i'}\cap f^{-1}(f(x))}$ at $x$ is a $C^r$ submersion onto the germ of $X_i\cap f^{-1}(\pi_j^Y\circ f(x))$ at $\pi_i^X(x)$. \par
  This definition of controlledness is stronger than that in \cite{G}. 
In \cite{G}, (4) is not assumed. 
However, if $f:\{X_i\}\to\{Y_j\}$ is a Thom map then (4) immediately follows from (1), (2) and (3), and existence of a $C^r$ tube system $\{T_i^X\}$ for $\{X_i\}$ controlled over a given controlled $C^r$ tube system $\{T_j^Y\}$ for $\{Y_j\}$ is known (see \cite{G} and \cite{S2}). 
We shall treat a $C^1$ stratification $f:\{X_i\}\to\{Y_j\}$ of $f$ which is not necessarily a Thom $C^1$ stratification but admits a controlled $C^1$ tube system $\{T_j^Y\}$ for $\{Y_j\}$ and a $C^1$ weak tube system $\{T_i^X\}$ for $\{X_i\}$ controlled over $\{T_j^Y\}$. \par
  In \cite{S3} theorem 1.1 is proved in the following more general form. 

\begin{thm}
Let $f:\{X_i\}\to\{Y_j\}$ be a $C^\infty$ stratification of a $C^\infty$ proper map $f:X\to Y$ between closed subsets of Euclidean spaces. 
Assume there exist a controlled $C^\infty$ tube system $\{T_j^Y\}$ for $\{Y_j\}$ and a $C^\infty$ tube system $\{T_i^X\}$ for $\{X_i\}$ controlled over $\{T_j^Y\}$. 
Then there exist homeomorphisms $\tau$ and $\pi$ from $X$ and $Y$ to polyhedra $P$ and $Q$, respectively, closed in some Euclidean spaces such that $\pi\circ f\circ\tau^{-1}:P\to Q$ is piecewise linear and $\tau(\overline{X_i})$ and $\pi(\overline{Y_j})$ are all polyhedra. 
If $f:\{X_i\}\to\{Y_j\},\,\{T_i^X\}$ and $\{T_j^Y\}$ are semialgebraic or, more generally, definable in an o-minimal structure, then we can choose semialgebraic or definable $\tau,\,\pi,\,P$ and $Q$. 
\end{thm}

(Note a semialgebraic closed polyhedron in a Euclidean space is semilinear, i.e., is defined by a finite number of equalities and inequalities of linear functions.) 
Moreover, the proof in \cite{S3} shows the following generalization though we do not repeat its proof. 

\begin{thm}
Let $f:\{X_i\}\to\{Y_j\}$ be a $C^1$ stratification of a $C^1$ proper map $f:X\to Y$ between closed subsets of Euclidean spaces. 
Let $I$ denote the set of indexes $i$ of $X_i$ such that $f|_{X_i}$ is not injective. 
Assume there exist a controlled $C^1$ tube system $\{T_j^Y\}$ for $\{Y_j\}$ and a $C^1$ weak tube system $\{T_i^X\}$ for $\{X_i\}$ controlled over $\{T_j^Y\}$ such that $\{T_i^X:i\in I\}$ is a $C^1$ tube system for $\{X_i:i\in I\}$. 
Then the result in theorem 2.1 holds. 
\end{thm}

We will prove theorem 1.2 by compactifying $f:X\to Y$ in theorem 1.2 and applying theorem 2.2 to the compctification. 
There are two unusual problems which we encounter. 
First the arguments do not work in the $C^2$ category and apply the $C^1$ category. 
Secondly we construct $\{T_j^Y:Y_j\subset\overline Y\}$ and $\{T_i^X:X_i\subset X\}$ by induction on $\dim Y_j$ and $\dim X_i$ but the induction of construction of $\{T_i^X:X_i\subset\overline X-X\}$ is downward. 
The two inductions are not independent and we need special conditions (iv) and (ix) for tube systems in the proof below. 
It is natural to ask whether we can extend $f$ to a Thom map $\overline f$. 
The answer is negative. 
To keep the property that $f$ is a Thom map also we use (iv) and (ix). 

\section{Proof theorem \ref{thmA}}
{\it Proof of theorem \ref{thmA}.} 
We assume $X$ is non-compact and $X$ and $Y$ are bounded in $\R^m$ and $\R^n$, respectively, by replacing $\R^m$ and $\R^n$ with $(0,\,1)^m$ and $(0,\,1)^n$ respectively. 
Then $\overline X-X$ and $\overline Y-Y$ are compact. 
Let $f:\{X_i\}\to\{Y_j\}$ be a semialgebraic Thom $C^1$ stratification of $f:X\to Y$. 
Then we can assume $f$ is extendable to $\overline X$. 
Apply Theorem II.4.1, \cite{S1} to the function on $\R^m$ measuring distance from the compact set $\overline X-X$. 
Then we have a non-negative semialgebraic $C^0$ function $\phi$ on $\R^m$ such that $\phi^{-1}(0)=\overline X-X$ and $\phi|_{\R^m-(\overline X-X)}$ is of class $C^1$. 
Choose $\epsilon>0\in\R$ so that $\phi$ is $C^1$ regular on $\phi^{-1}((0,\,\epsilon])$ and let $\phi'$ be a semialgebraic $C^1$ function on $\R$ such that $\phi'(0)=0$, $\phi'$ is regular on $(0,\,\epsilon)$ and $\phi'=1$ on $[\epsilon,\,\infty)$. 
Set 
$$
\Phi(x)=(\phi'\circ\phi(x),\phi'\circ\phi(x)x)\quad\text{for}\ x\in X. $$
Then $\Phi$ is a semialgebraic $C^1$ imbedding of $X$ into $\R^{m+1}$ such that $\Phi(X)$ is bounded and $\overline{\Phi(X)}-\Phi(X)=\{0\}$. 
Hence replacing $X$ with $\Phi(X)$ we assume $\overline X-X=\{0\}$ from the beginning. 
Moreover, replace $X$ with the graph of $f$.
Then we suppose $X$ is contained and bounded in $\R^m\times\R^n$, $\overline X-X\subset\{0\}\times\overline Y$, $f:X\to Y$ is the restriction of the projection $p:\R^m\times\R^n\to\R^n$ and hence $f$ is extended to a semialgebraic $C^1$ map $\overline f:\overline X\to\overline Y$. \par
  By the same reason we assume $\overline Y-\{0\}$. 
Note then $\{Y_j,0\}$ is a semialgebraic Whitney $C^1$ stratification of $\overline Y$. 
Let $\{T_j^Y\}$ be a controlled semialgebraic $C^1$ tube system for $\{Y_j\}$ and $\{T_i^X\}$ a semialgebraic $C^1$ tube system for $\{X_i\}$ controlled over $\{T_j^Y\}$. 
Assume the set of indexes of $Y_j$ does not contain 0, set $Y_0=\{0\}$ and add $Y_0$ to $\{Y_j\}$. 
Then we can assume there is a semialgebraic $C^1$ tube $T_0^Y=(|T_0^Y|,\pi_0^Y,\rho_0^Y)$ at $Y_0$ such that $\{T_j^Y,T_0^Y:j\not=0\}$ is controlled for the following reason. \par
  Let $|T_0^Y|$ be the closed ball $B(\epsilon)$ with center 0 in $\R^n$ and with small radius $\epsilon>0$ (we treat closed balls in place of open balls for simplicity of notation), and set $\pi_0^Y(y)=0$ and, tentatively, $\rho_0^Y(y)=|y|^2$ for $y\in|T_0^Y|$. 
Then the condition $\rho_0^Y\circ\pi_j^Y=\rho_0^Y$ on $|T_0^Y|\cap|T_j^Y|$ for $j\not=0$ does not necessarily hold. 
For that condition it suffices to find a semialgebraic homeomorphism $\tau$ of $\R^n$ of class $C^1$ outside of 0 and such that $\tau(0)=0$, $\tau=\id$ outside of $B(\epsilon)$ and $\rho_0^Y\circ\pi_j^Y\circ\tau^{-1}=\rho_0^Y$ on $B(\epsilon')\cap\tau(|T_j^Y|)$ for $j\not=0$, shrunk $|T_j^Y|$ and some $\epsilon'>0$. \par
  Let $Y_j$ be such that $\dim Y_j$ is the smallest in $\{Y_j:0\in\overline{Y_j},\,j\not=0\}$, and choose $\epsilon$ so small that $\rho_0^Y|_{Y_j\cap|T_0^Y|}$ is $C^1$ regular, which implies that $\rho_0^{Y-1}(\epsilon'^2)$ is transversal to $Y_j$ for any $0<\epsilon'\le\epsilon$. 
Set $Y_j(\epsilon')=Y_j\cap\rho_0^{Y-1}(\epsilon'^2)$. 
We will define a semialgebraic homeomorphism $\tau_j$ of $\R^n$ of class $C^1$ outside of 0 such that $\tau_j(0)=0$, $\tau_j=\id$ outside of $B(\epsilon)$ and $\rho_0^Y\circ\pi_j^Y\circ\tau_j^{-1}=\rho_0^Y$ on $B(\epsilon/2)\cap\tau_j(|T_j^Y|)$ for shrunk $|T_j^Y|$. 
Since the problem is local at $Y_j$, we can assume by Thom's first isotopy lemma (see Theorem II.6.1 and it complement, [4]) that 
$$
|T_0^Y|\cap Y_j=Y_j(\epsilon)\times(0,\,\epsilon^2],\ \text{after then,}\ |T_0^Y|\cap|T_j^Y|=\cup\{y+L_y:y\in Y_j(\epsilon)\}\times(0,\,\epsilon^2]
$$
and $\pi_j^Y(y+z,t)$ and $\rho_0^Y(y+z,t)$ are of the form $(y,\pi_j^{Y\prime}(y+z,t))$ and $t$, respectively, for $y\in Y_j(\epsilon)$ and $(z,t)\in L_y\times(0,\,\epsilon^2]$, where $L_y$ is a linear subspace of the tangent space $T_y\rho_0^{Y-1}(\epsilon^2)$ of codimension\,=\,$\codim Y_j$ in $\R^n$ such that the correspondence $Y_j(\epsilon)\ni y\to L_y\in G_{n,\codim Y_j}$ is semialgebraic and of class $C^1$ and $\pi_j^{Y\prime}$ is a semialgebraic $C^1$ function defined on $\cup\{y+L_y\}\times(0,\,\epsilon^2]$. 
For simplicity of notation we write $\cup_{y\in Y_j(\epsilon)}\{y\}\times L_y$ as $Y_j(\epsilon)\times L$. 
Transform $Y_j(\epsilon)\times L\times(0,\,\epsilon^2]$ by a semialgebraic $C^1$ diffeomorphism $(y,z,t)\to(y,z/kt^k,t)$ for sufficiently large integer $k$. 
Then we can assume 
$$|\pi_j^{Y\prime}(y+z,t)-t|\le\epsilon^2/28\ \ \text{and}\ \ |\frac{\partial\pi_j^{Y\prime}}{\partial t}(y+z,t)-1|<1/4\ \ \text{for}\ |z|\le 1\leqno{(0)}
$$
$$
\pi_j^{Y\prime}(y,t)=t. \leqno{\text{since}}
$$
Let $\xi$ be a semialgebraic $C^1$ function on $\R$ such that $0\le\xi\le 1,\ \xi=1$ on $(-\infty,\,1/2),\ \xi=0$ on $(2/3,\,\infty)$ and $|\frac{d\xi}{d t}|\le 7$. 
Set
$$
\tau_j(y+z,t)=(y+z,(1-\xi(2t/\epsilon^2)\xi(|z|))t+\xi(2t/\epsilon^2)\xi(|z|)\pi_j^{Y\prime}(y+z,t))\qquad
$$
$$
\qquad\qquad\qquad\qquad\qquad\qquad\qquad\qquad\qquad\text{for}\ (y,z,t)\in Y_j(\epsilon)\times L\times(0,\,\epsilon^2]. 
$$
Then $\tau_j=\pi_j^Y$ if $t\le\epsilon^2/4$ and $|z|\le 1/2$, $\tau_j=\id$ if $t\ge\epsilon^2/3$ or $|z|\ge2/3$ and, moreover, $\tau_j$ is a diffeomorphism because 
$$
|\frac{\partial}{\partial t}\big((1-\xi(t/\epsilon^2)\xi(|z|))t+\xi(t/\epsilon^2)\xi(|z|)\pi_j^{Y\prime}(y+z,t)\big)-1|\qquad\qquad\qquad\quad
$$
$$
\le\xi(t/\epsilon^2)\xi(|z|)|1-\frac{\partial\pi_j^{Y\prime}}{\partial t}(y+z,t)|+|\frac{d\xi}{d t}(t/\epsilon^2)\xi(|z|)|t-\pi_j^{Y\prime}(y+z,t)|/\epsilon^2
$$
$$
\qquad\qquad\qquad\qquad\qquad\qquad\le 1/4+1/4=1/2\quad\text{for}\ |z|\le 1. 
$$
Thus we can assume $\rho_0^Y\circ\pi_j^Y=\rho_0^Y$ on $|T_0^Y|\cap|T_j^Y|$. \par
  Repeating the same arguments by induction on $\dim Y_{j'}$ for all $Y_{j'}$ with $0\in\overline{Y_{j'}}$ we obtain the required $\tau$. 
Here we note only that for $j'$ with $\overline{Y_{j'}}-Y_{j'}\supset Y_j$, though $Y_{j'}(\epsilon)$ is not compact, $(0)$ can holds. 
Indeed 
$$
\rho_0^Y=\rho_0^Y\circ\pi_j^Y\circ\pi_{j'}^Y=\rho_0^Y\circ\pi_{j'}^Y\quad\text{on}\ |T_0^Y|\cap|T_j^Y|\cap|T_{j'}^Y|. 
$$
Hence when we describe $\pi_{j'}^Y$ as above there is a semialgebraic neighborhood $U$ of $Y_j(\epsilon)\times(0,\,\epsilon^2]$ in $\overline{Y_{j'}(\epsilon)}\times(0,\,\epsilon^2]$ such that 
$$
\pi_{j'}^{Y\prime}(y+z,t)=t\quad\text{for}\ (y,z,t)\in Y_{j'}(\epsilon)\times L_y\times(0,\,\epsilon^2]\ \text{with}\ (y,t)\in U. 
$$
In conclusion we assume $Y$ is compact. \par
  If $f:\{X_i\}\to\{Y_j\}$ is extended to a Thom $C^1$ stratification of $\overline f:\overline X\to Y$, then theorem 1.2 follows from theorem 1.1 in the $C^1$ case. 
However, such extension does not always exist. 
Instead we will find a semialgebraic $C^1$ stratification $\overline f:\{X'_{i'}\}\to\{Y'_{j'}\}$ of $\overline f$, a controlled semialgebraic $C^1$ tube system $\{T_{j'}^{Y\prime}\}$ for $\{Y'_{j'}\}$ and a semialgebraic $C^1$ weak tube system $\{T_{i'}^{X\prime}\}$ for $\{X'_{i'}\}$ controlled over $\{T_{j'}^{Y\prime}\}$ such that $\{X'_{i'}\}|_X$ and $\{Y'_{j'}\}|_Y$ are substratifications of $\{X_i\}$ and $\{Y_j\}$. 
Here $\{Y'_{j'}\}$ is a Whitney stratification but $\{X'_{i'}\}$ is not necessarily so. \par
  Set $Z=\overline X-X$, which is compact. 
Note $Z=\{0\}\times\overline f(Z)$ and $\overline f|_Z$ is a homeomorphism onto $\overline f(Z)$. 
Let $\{Y'_{j'}\}$ be a semialgebraic Whitney $C^1$ substratification of $\{Y_j\}$ such that each stratum is connected, $\overline f(Z)$ is a union of some $Y'_{j'}$'s and $\{X_i,\{0\}\times(Y'_{j'}\cap\overline f(Z))\}$ is a Whitney $C^1$ stratification of $\overline X$, which is constructed in the same way as the canonical semialgebraic $C^\omega$ stratification of a semialgebraic set since $\overline f(Z)$ is closed in $Y$. 
Note $\{Y'_{j'}\}$ satisfies the frontier condition. 
Set 
$$
\{X'_{i'}\}=\{X_i\cap\overline f^{-1}(Y'_{j'}),Z\cap\{0\}\times Y'_{j'}\}. 
$$
Then $\{X'_{i'}\}$ is a semialgebraic (not necessarily Whitney) $C^1$ stratification of $\overline X$; $\{X'_{i'}\cap X\}$ is a substratification of $\{X_i\}$; $\overline f:\{X'_{i'}\}\to\{Y'_{j'}\}$ is a $C^1$ stratification of $\overline f$; we can choose $\{Y'_{j'}\}$ so that for each $Y'_{j'}$, $\{X'_{i'}:\overline f(X'_{i'})=Y'_{j'}\}$ is a Whitney $C^1$ stratification for the following reason. \par
  Assume $Y'_{j'}\not\subset\overline f(Z)$. 
Then $Y'_{j'}\cap\overline f(Z)=\emptyset$ and there is $Y_j$ including $Y'_{j'}$. 
By definition of $\{X'_{i'}\}$, 
$$
\{X'_{i'}:\overline f(X'_{i'})=Y'_{j'}\}=\{X_i\cap f^{-1}(Y'_{j'})\}. 
$$
Therefore the assertion follows from the fact that given a Whitney $C^r$ stratification $\{M_1,M_2\}$, a $C^r$ map $g$ from $M_1\cup M_2$ to a $C^r$ manifold $N$ such that $g|_{M_1}$ and $g|_{M_2}$ are $C^r$ submersions into $N$ and a $C^r$ submanifold $N_1$ of $N$ then $\{M_1\cap g^{-1}(N_1),M_2\cap g^{-1}(N_1)\}$ is a Whitney $C^r$ stratification. \par
  Next assume $Y'_{j'}\subset\overline f(Z)$, and let $X'_{i'_1}$ and $X'_{i'_2}$ be such that $\overline f(X'_{i'_k})=Y'_{j'},\ k=1,2$, and $(\overline{X'_{i'_1}}-X'_{i'_1})\cap X'_{i'_2}\not=\emptyset$. 
Then we need to see $(X'_{i'_1},X'_{i'_2})$ can satisfy the Whitney condition. 
Since $\overline f|_Z$ is injective, there are only two possible cases to consider: $X'_{i'_k}=X_{i_k}\cap\overline f^{-1}(Y'_{j'}),\ k=1,2$, for some $i_1$ and $i_2$ or $X'_{i'_1}=X_{i_1}\cap\overline f^{-1}(Y'_{j'})$ and $X'_{i'_2}=\{0\}\times Y'_{j'}$. 
In the former case there is $j$ such that $Y'_{j'}\subset Y_j$. 
Hence the Whitney condition is satisfied by the same reason as in the case of $Y'_{j'}\not\subset\overline f(Z)$. 
Consider the latter case. 
If $\{X'_{i'_1},\{0\}\times Y'_{j'}\}$ is not a Whitney stratification, let $Y''_{j'}$ denote the subset of $Y'_{j'}$ consisting of $y$ such that $(X'_{i_1},\{0\}\times Y'_{j'})$ does not satisfy the Whitney condition at $(0,y)$. 
Then $Y''_{j'}$ and hence $\overline{Y''_{j'}}$ are semialgebraic and of dimension smaller that $\dim Y'_{j'}$. 
Divide $Y'_{j'}$ to $\{Y'_{j'}-\overline{Y''_{j'}},\overline{Y''_{j'}}\}$ and substratify $\{Y'_{j'}\cap\overline f(Z)\}$ by downward induction on dimension of $Y'_{j'}$ so that the above conditions on $\{Y'_{j'}\}$ are kept and $Y''_{j'}=\emptyset$. 
Then $\{X'_{i'_1},\{0\}\times Y'_{j'}\}$ becomes a Whitney stratification. \par
  Now we define a controlled semialgebraic $C^1$ tube system $\{T_{j'}^{Y\prime}=(|T_{j'}^{Y\prime}|,\pi_{j'}^{Y\prime},\rho_{j'}^{Y\prime})\}$ for $\{Y'_{j'}\}$. 
For simplicity of notation, assume $\dim Y_j=j$ gathering strata of the same dimension. 
For each $j$, set 
$$
J_j=\left\{
\begin{array}{ll}
\{j':Y'_{j'}\subset Y_j,\}\quad&\text{if}\ j\ge0,\\
\emptyset\quad&\text{if}\ j=-1. \\
\end{array}
\right.
$$
We define $\{T_{j'}^{Y\prime}:j'\in J_j\}$ by induction on $j$. 
Fix a non-negative integer $j_0$, and assume we have constructed a controlled semialgebraic $C^1$ tube system $\{T_{j'}^{Y\prime}:j'\in J_j,j<j_0\}$ so that $T_{j'}^{Y\prime}=T_{j_1}^Y|_{|T_{j'}^{Y\prime}|}$ for $j'\in J_{j_1},\ j_1<j_0$, with $\dim Y'_{j'}=j_1$,
$$
\pi_{j'}^{Y\prime}\circ\pi_j^Y=\pi_{j'}^{Y\prime}\quad\text{on}\ |T_{j'}^{Y\prime}|\cap|T_j^Y|\ \text{for}\ j'\ \text{and}\ j\ \text{with}\ Y'_{j'}\subset\overline{Y_j},\leqno{(*)_Y}
$$
$$
\rho_{j'}^{Y\prime}\circ\pi_j^Y=\rho_{j'}^{Y\prime}\quad\text{on}\ |T_{j'}^{Y\prime}|\cap|T_j^Y|\ \text{for}\ j'\in J_{j_1}\ \text{and}\ j\ \text{with}\ j_1<j,\leqno{(**)_Y}
$$
$\pi_{j'}^{Y\prime}$ are of class $C^1$ and $\rho_{j'}^{Y\prime}$ are of class $C^1$ on $|T_{j'}^{Y\prime}|-Y'_{j'}$. 
For the conditions of the first and $(**)_Y$ we need to proceed in the $C^1$ category because there does not necessarily exist such $\{T_{j'}^{Y\prime}\}$ of class $C^2$ even if $\{T_j^Y\}$ is of class $C^2$. \par
  We wil define a semialgebraic $C^1$ tube system $\{T_{j'}^{Y\prime}:j'\in J_{j_0}\}$ for $\{Y'_{j'}:j'\in J_{j_0}\}$. 
For the time being, let $\{T_{j'}^{Y\prime}:j'\in J_{j_0}\}$ be a semialgebraic $C^1$ tube system for $\{Y'_{j'}:j'\in J_{j_0}\}$ such that $\{T_{j'}^{Y\prime}:j'\in J_j,j\le j_0\}$ is controlled (Lemma II.6.10, \cite{S2} states only the case where $\cup_{j'\in J_{j_0}}Y'_{j'}$ is compact but its proof works in the general case. 
We omit the details.) 
We modify $\{T_{j'}^{Y\prime}:j'\in J_{j_0}\}$ so that the conditions are satisfied. 
Let $j'\in J_{j_0}$. \par
  Restrict $\pi_{j'}^{Y\prime}$ and $\rho_{j'}^{Y\prime}$ to $Y_{j_0}$ for $j'\in J_{j_0}$ and define afresh them outside of $Y_{j_0}$ as follows. 
Let $\pi_{j'}^{Y\prime}$ and $\rho_{j'}^{Y\prime}$, $j'\in J_{j_0}$, now denote the restrictions. 
If $\dim Y'_{j'}=j_0$, we should set $T_{j'}^{Y\prime}=T_{j_0}^Y|_{|T_{j'}^{Y\prime}|}$. 
Then $(*)_Y$ and $(**)_Y$ are satisfied because $\{Y_j^Y\}$ is controlled. 
Assume $\dim Y'_{j'}<j_0$ and hence $j_0>0$. 
In this case, define the extension of $\pi_{j'}^{Y\prime}$ to $|T_{j'}^{Y\prime}|$ to be $\pi_{j'}^{Y\prime}\circ\pi_{j_0}^Y$, and keep the same notation $\pi_{j'}^{Y\prime}$ for the extension. 
Then by controlledness of $\{T_j^Y\}$, $(*)_Y$ holds for any $j$ with $Y'_{j'}\subset\overline{Y_J}$. 
The problem is how to extend $\rho_{j'}^{Y\prime}$. \par
  As the problem is local at $Y'_{j'}$ (see II.1.1, [4]), considering semialgebraic tubular neighborhoods of $Y'_{j'}$ and $Y_{j_0}$ we can assume for each $y\in Y'_{j'}$, $\pi_{j'}^{Y\prime-1}(y),\ \pi_{j'}^{Y\prime-1}(y)\cap Y_{j_0}$ and $\pi_{j_0}^{Y-1}(y)$ are of the form $y+L_y,\ y+L_{0,y}$ and $y+L_{0,y}^\perp$, where $L_y$ and $L_{0,y}$ are linear subspaces of $\R^n$ with $L_y\supset L_{0,y}$ and $L_{0,y}^\perp$ is the orthocomplement of $L_{0,y}$ with respect to $L_y$, and $\pi_{j_0}^Y|_{\pi_{j'}^{Y\prime-1}(y)}:\pi_{j'}^{Y\prime-1}(y)\longrightarrow\pi_{j'}^{Y\prime-1}(y)\cap Y_{j_0}$ is induced by the orthogonal projection of $L_y$ to $L_{0,y}$ and 
$$
  \rho_{j_0}^Y(y+z_1+z_2)=|z_2|^2\ \text{for}\ (y,z_1,z_2)\in Y'_{j'}\times L_{0,y}\times L_{0,y}^\perp, 
$$
where $Y'_{j'}\times L_{0,y}\times L_{0,y}^\perp$ denotes $\cup_{y\in Y'_{j'}}\{y\}\times L_{0,y}\times L_{0,y}^\perp$. 
$$
\rho_{j'}^{Y\prime\prime}(y+z_1+z_2)=|z_1|^2+|z_2|^2\quad\text{for}\ (y,z_1,z_2)\in Y'_{j'}\times L_{0,y}\times L_{0,y}^\perp. \leqno{\text{Set}}
$$
Then $(|T_{j'}^{Y\prime}|,\pi_{j'}^{Y\prime},\rho_{j'}^{Y\prime\prime})$ is a semialgebraic $C^1$ tube at $Y'_{j'}$ but not always satisfy the condition $\rho_{j'}^{Y\prime\prime}\circ\pi_{j_0}^Y=\rho_{j'}^{Y\prime\prime}$. 
We need to modify $\rho_{j'}^{Y\prime\prime}$ so that the equality holds on a neighborhood of $Y_{j_0}-Y'_{j'}$. 
Let $\xi$ be a semialgebraic $C^1$ function on $\R$ such that $\xi=1$ on $(-\infty,\,1],\ \xi=0$ on $[2,\infty)$ and $d\xi/d t\le0$. 
Set 
$$
\eta_{j'}(z_1,z_2)=\left\{
\begin{array}{ll}
\xi(\frac{|z_2|}{|z_1|^2})\frac{|z_1|}{(|z_1|^2+|z_2|^2)^{1/2}}+1-\xi(\frac{|z_2|}{|z_1|^2})\, &\text{for}\ (z_1,z_2)\in(L_{0,y}-\{0\})\times L_{0,y}^\perp,\\
1\ &\text{for}\ (z_1,z_2)\in\{0\}\times L_{0,y}^\perp,\\
\end{array}
\right.
$$
and define a semialgebraic map $\tau_{j'}$ between $|T_{j'}^{Y\prime}|$ by 
$$
\tau_{j'}(y+z_1+z_2)=y+\eta_{j'}(z_1,z_2)z_1+\eta_{j'}(z_1,z_2)z_2\quad\text{for}\ (y,z_1,z_2)\in Y'_{j'}\times L_{0,y}\times L_{0,y}^\perp. 
$$\par
  Then $\pi_{j'}^{Y\prime}\circ\tau_{j'}=\pi_{j'}^{Y\prime}$; 
$$
\tau_{j'}=\id\quad\text{on}\ \{y+z_1+z_2:|z_2|\ge2|z_1|^2\};
$$
$$
\tau_{j'}(y+z_1+z_2)=y+\frac{|z_1|}{(|z_1|^2+|z_2|^2)^{1/2}}z_1+\frac{|z_1|}{(|z_1|^2+|z_2|^2)^{1/2}}z_2\qquad\qquad\qquad
$$
$$\qquad\qquad\qquad\qquad\quad\text{for}\ (y,z_1,z_2)\in Y'_{j'}\times L_{0,y}\times L_{0,y}^\perp\ \text{with}\ |z_2|\le|z_1|^2;
$$
$$
\rho_{j'}^{Y\prime\prime}\circ\tau_{j'}(y+z_1+z_2)=|z_1|^2\quad\text{for the same}\ (y,z_1,z_2); \leqno{(***)_Y}
$$
for each line $l$ in $\{y\}\times L_{0,y}\times L_{0,y}^\perp$ passing through 0 parameterized by $t\in\R$ as $z_1=z_1(t)$ and $z_2=z_2(t)$ so that $|z_1(t)|=|t|$ and $|z_2(t)|=a|t|$ for $a\ge 0\in\R$, 
$$
\tau_{j'}(l)=l,
$$
$$
|\tau_{j'}(y+z_1(t)+z_2(t))-y|=\eta_{j'}(z_1(t),z_2(t))(|z_1(t)|^2+|z_2(t)|^2)^{1/2}\quad\qquad\qquad
$$
$$
\qquad\qquad\qquad\qquad\qquad\qquad\qquad\qquad=\xi(\frac{a}{|t|})|t|+(1-\xi(\frac{a}{|t|}))(1+a^2)^{1/2}|t|,
$$
hence by easy calculations we see if $a$ is sufficiently small then $\tau_{j'}|_l$ is a $C^1$ diffeomorphism of $l$ and, therefore by the above equality $\tau_{j'}=\id$ on $\{|z_2|\ge 2|z_1|^2\}$ shrinking $|T_{j'}^{Y\prime}|$ we can assume $\tau_{j'}$ is a homeomorphism and its restriction to $|T_{j'}^{Y\prime}|-Y'_{j'}$ is a $C^1$ diffeomorphism; 
moreover, if we set $\rho_{j'}^{Y\prime}=\rho_{j'}^{Y\prime\prime}\circ\tau_{j'}$ and $T_{j'}^{Y\prime}=(|T_{j'}^{Y\prime}|,\pi_{j'}^{Y\prime},\rho_{j'}^{Y\prime})$ for all $j'\in J_{j_0}$ with $\dim Y'_{j'}<j_0$ then $\{T_{j'_1}^{Y\prime}:j'_1\in J_{j_1},j_1\le j_0\}$ is a controlled semialgebraic $C^1$ tube system. 
Indeed, for $j'_1\in J_{j_0}$ and $j'_2$ with $(\overline{Y'_{j'_1}}-Y'_{j'_1})\cap Y'_{j'_2}\not=\emptyset$, the following equalities folds on $|T_{j'_1}^{Y\prime}|\cap|T_{j'_2}^{Y\prime}|$ 
$$
\pi_{j'_2}^{Y\prime}\circ\pi_{j'_1}^{Y\prime}=\pi_{j'_2}^{Y\prime}\circ\pi_{j'_1}^{Y\prime}\circ\pi_{j_0}^Y\quad\text{by definition of}\ \pi_{j'_1}^{Y\prime}\qquad\qquad\quad
$$
$$
=\pi_{j'_2}^{Y\prime}\circ\pi_{j_0}^Y\quad\text{by controlledness of}\ \{T_{j'}^{Y\prime}|_{Y_{j_0}}:j'\in J_j,\,j\le j_0\}
$$
$$
=\pi_{j'_2}^{Y\prime}\ \text{by definition of}\ \pi_{j'_2}^{Y\prime}\ \text{in the case of }j'_2\in J_{j_0}\text{ and by }(*)_Y\text{ in the other case}. 
$$
In the same way we see by $(**)_Y$ and $(***)_Y$ 
$$
\rho_{j'_2}^{Y\prime}\circ\pi_{j'_1}^{Y\prime}=\rho_{j'_2}^{Y\prime}\quad\text{on}\ |T_{j'_1}^{Y\prime}|\cap|T_{j'_2}^{Y\prime}|. 
$$
Hence it remains to show $\tau_{j'}$ is a $C^1$ diffeomorphism. \par
  It is easy to show $\tau_{j'}$ is differentiable at $Y'_{j'}$ and its differential $d\tau_{j'a}$ at each point $a$ of $Y'_{j'}$ is equal to the identity map. 
Hence we only need to show the map $|T_{j'}^{Y\prime}|\ni a\to d\tau_{j'a}\in GL(\R^n)$ is of class $C^0$. 
As the problem is local at each point of $Y'_{j'}$ we suppose 
$$
Y'_{j'}=\R^{n'}\times\{0\}\times\{0\},\ Y_{j_0}=\R^{n'}\times\R^{n_1}\times\{0\},\ |T_{j'}^{Y\prime}|=|T_{j_0}^Y|=\R^{n'}\times\R^{n_1}\times\R^{n_2}
$$
and $\pi_{j_0}^Y$ and $\pi_{j'}^{Y\prime}$ are the projections of $\R^{n'}\times\R^{n_1}\times\R^{n_2}$ to $\R^{n'}\times\R^{n_1}\times\{0\}$ and $\R^{n'}\times\{0\}\times\{0\}$ respectively. 
Then it suffices to see the differential at $(z_{01},z_{02})$ of the map $\R^{n_1}\times\R^{n_2}\ni(z_1,z_2)\to(\eta_{j'}(z_1,z_2)z_1,\eta_{j'}(z_1,z_2)z_2)\in\R^{n_1}\times\R^{n_2}$ converges to the identity map as $(z_{01},z_{02})\to(0,0)$. 
That is, 
$$
\aligned
&d\left(\frac{\xi(\frac{|z_2|}{|z_1|^2})((|z_1|^2+|z_2|^2)^{1/2}-|z_1|)z_i}{(|z_1|^2+|z_2|^2)^{1/2}}\right)_{(z_{01},z_{02})}=\\
&d\left(\frac{\xi(\frac{|z_2|}{|z_1|^2})|z_2|^2z_i}{(|z_1|^2+|z_2|^2)^{1/2}((|z_1|^2+|z_2|^2)^{1/2}+|z_1|)}\right)_{(z_{01},z_{02})}\longrightarrow 0
\endaligned
$$
as $(z_{01},z_{02})\to(0,0)$ with $|z_2|\le2|z_1|^2,\ i=1,2,$ since $\eta_{j'}(z_1,z_2)=1$ for $(z_1,z_2)$ with $|z_2|\ge2|z_1|^2$. 
That is easy to check. 
We omit the details. \par
  Thus we obtain semialgebraic $C^1$ tubes $T_{j'}^{Y\prime}$ for all $j'\in J_{j_0}$. 
The other requirements in the induction hypothesis are satisfied as follows. 
By definition of $T_{j'}^{Y\prime}$, 
$$
T_{j'}^{Y\prime}=T_{j_0}^Y|_{|T_{j'}^{Y\prime}|}\quad\text{for}\ j'\in J_{j_0}\ \text{with}\ \dim Y'_{j'}=j_0;
$$ by controlledness of $\{T_j^Y\}$ and by definition of $T_{j'}^{Y\prime}$, for $j'$ and $j$ with $Y'_{j'}\subset\overline{Y_j},\ j'\in J_{j_0}$ and $j\ge j_0$, 
$$
\pi_{j'}^{Y\prime}\circ\pi_j^Y=\pi_{j'}^{Y\prime}\circ\pi_{j_0}^Y\circ\pi_j^Y=\pi_{j'}^{Y\prime}\circ\pi_{j_0}^Y=\pi_{j'}^{Y\prime}\quad\text{on}\ |T_{j'}^{Y\prime}|\cap|T_j^Y|;\leqno{(*)_Y}
$$
$(**)_Y$ holds for $j'$ and $j$ with $j'\in J_{j_0}$ and $j>j_0$ for the following reason. \par
  That is clear if $\dim Y'_{j'}=j_0$. 
Hence assume $\dim Y'_{j'}<j_0$ and use the above coordinate system $Y'_{j'}\times L_{0,y}\times L_{0,y}^\perp$. 
Then 
$$
\rho_{j'}^{Y\prime}(y+z_1+z_2)\!=\!\rho_{j'}^{Y\prime\prime}\circ\tau_{j'}(y+z_1+z_2)\!=\!\eta_{j'}^2(z_1,z_2)(|z_1|^2+|z_2|^2)\qquad
$$
$$
\qquad\qquad\qquad\qquad\qquad\qquad\qquad\qquad\text{for }(y,z_1,z_2)\in Y'_{j'}\times L_{0,y}\times L_{0,y}^\perp 
$$
and $\eta_{j'}(z_1,z_2)$ depends on only $|z_1|$ and $|z_2|$. 
Hence if we set 
$$
\pi_j^Y(y+z_1+z_2)=\pi_{j1}^Y(y+z_1+z_2)+\pi_{j2}^Y(y+z_1+z_2)+\pi_{j3}^Y(y+z_1+z_2),
$$
$$
\pi_{j1}^Y(y+z_1+z_2)\in Y'_{j'},\ \pi_{j2}^Y(y+z_1+z_2)\in L_{0,y},\ \pi_{j3}^Y(y+z_1+z_2)\in L_{0,y}^\perp. 
$$
then it suffices to see 
$$
\pi_{j2}^Y(y+z_1+z_2)=z_1\quad\text{and}\quad|\pi_{j3}^Y(y+z_1+z_2)|=|z_2|. 
$$
By controlledness of $\{T_j^Y\}$ we have $\pi_{j_0}^Y\circ\pi_j^Y=\pi_{j_0}^Y$. 
Hence by the equation $\pi_{j_0}^Y(y+z_1+z_2)=y+z_1$, the former equality holds. 
The latter also follows from the equations $\rho_{j_0}^Y\circ\pi_j^Y=\rho_{j_0}^Y$ and $\rho_{j_0}^Y(y+z_1+z_2)=|z_2|^2$. \par
  Hence by induction we have a controlled semialgebraic $C^1$ tube system $\{T_{j'}^{Y\prime}\}$ for $\{Y'_{j'}\}$ such that $T_{j'}^{Y\prime}=T_j^Y|_{|T_{j'}^{Y\prime}|}$ for $j'\in J_j$ with $\dim Y'_{j'}=j$, $(*)_Y$ for $j'$ and $j$ with $Y'_{j'}\subset\overline{Y_j}$ and $(**)_Y$ for $j'\in J_{j_1}$ and $j$ with $j_1<j$. \par
  Next we define $\{T_{i'}^{X\prime}\}$ by induction as $\{T_{j'}^{Y\prime}\}$. 
Consider all $X'_{i'}$ included in $X$ and forget $X'_{i'}$ outside of $X$. 
We change the set of indexes of $X_i$. 
For non-negative integers $i_0$ and $j_0$, let $X_{i_0,j_0}$ denote the union of $X_i$'s such that $\dim X_i=i_0$ and $f(X_i)\subset Y_{j_0}$, i.e., $\dim f(X_i)=j_0$, naturally define $T_{i,j}^X=(|T_{i,j}^X|,\pi_{i,j}^X,\rho_{i,j}^X)$ and continue to define $\{X'_{i'}\}$ to be $\{X_{i,j}\cap p^{-1}(Y'_{j'}),Z\cap\{0\}\times Y'_{j'}\}$. 
Then $\dim X_{i,j}=i$ and $f|_{X_{i,j}}$ is a map to $Y_j$. 
Let $I_i$ denote the set of indexes of $X'_{i'}$ such that $X'_{i'}$ is included in $X_{i,j}$ for some $j$. 
Note $X=\cup\{X'_{i'}:i'\in I_i$ for some $i\}$. 
Fix a non-negative integer $i_0$, and assume there exists a semialgebraic $C^1$ tube system $\{T_{i'}^{X\prime}=(|T_{i'}^{X\prime}|,\pi_{i'}^{X\prime},\rho_{i'}^{X\prime}):i'\in I_i,\ i<i_0\}$ for $\{X'_{i'}:i'\in I_i,\ i<i_0\}$ such that the following four conditions are satisfied, which are, except (iv), similar to the conditions (1), (2) and (3) in section 2. \par\noindent
(i) For $i,\ i'$ and $j'$ with $i<i_0,\ i'\in I_i$ and $f(X'_{i'})=Y'_{j'}$, 
$$
f\circ\pi_{i'}^{X\prime}=\pi_{j'}^{Y\prime}\circ p\quad\text{on}\ |T_{i'}^{X\prime}|\cap p^{-1}(|T_{j'}^{Y\prime}|). 
$$
(ii) For each $j'$, $\{T_{i'}^{X\prime}:f(X'_{i'})=Y'_{j'},\ i'\in I_i,\ i<i_0\}$ is a controlled semialgebraic $C^1$ tube system for $\{X'_{i'}:f(X'_{i'})=Y'_{j'},\ i'\in I_i,\ i<i_0\}$. \par\noindent
(iii) For $i_k,\,i'_k,\ k=1,2,3,\ i_4$ and $j_4$ with $i_k<i_0,\ i'_k\in I_{i_k},\ k=1,2,3$, $X'_{i'_1}\cap(\overline{X'_{i'_2}}-X'_{i'_2})\not=\emptyset$ and $X'_{i'_3}\subset\overline{X_{i_4,j_4}}$, 
$$
\pi_{i'_1}^{X\prime}\circ\pi_{i'_2}^{X\prime}=\pi_{i'_1}^{X\prime}\quad\text{on}\ |T_{i'_1}^{X\prime}|\cap|T_{i'_2}^{X\prime}|, 
$$
$$
\pi_{i'_3}^{X\prime}\circ\pi_{i_4,j_4}^X=\pi_{i'_3}^{X\prime}\quad\text{on}\ |T_{i'_3}^{X\prime}|\cap|T_{i_4,j_4}^X|, 
$$
if $i_3<i_4$ moreover, then
$$
\rho_{i'_3}^{X\prime}\circ\pi_{i_4,j_4}^X=\rho_{i'_3}^{X\prime}\quad\text{on}\ |T_{i'_3}^{X\prime}|\cap|T_{i_4,j_4}^X|. 
$$
(iv) For $i,\ i'$ and $j$ with $i<i_0,\ i'\in I_i$ and $\dim X'_{i'}=i$, 
$$
T_{i'}^{X\prime}=T_{i,j}^X|_{|T_{i'}^{X\prime}|}. 
$$\par
  Then we need to define $\{T_{i'}^{X\prime}:i'\in I_{i_0}\}$ so that the induction process works. 
Before that we note a fact. \par\noindent
(v) Given $i_k,\,i'_k,\,j'_k,\ k=1,2$, with $i_k<i_0,\ i'_k\in I_{i_k},\ k=1,2,\ X'_{i'_1}\cap(\overline{X'_{i'_2}}-X'_{i'_2})\not=\emptyset$, $Y'_{j'_1}\subset\overline{Y'_{j'_2}}-Y'_{j'_2}$ and $f(X'_{i'_k})=Y'_{j'_k},\ k=1,2$, then the restriction of the map $(\pi_{i'_1}^{X\prime},f)$ to $X'_{i'_2}\cap|T_{i'_1}^{X\prime}|$ is a $C^1$ submersion into the fiber product $X'_{i'_1}\times_{(f,\pi_{j'_1}^{Y\prime})}(Y'_{j'_2}\cap|T_{j'_1}^{Y\prime}|)$. \par
The reason is the following. \par
  Case where $X'_{i'_k}\subset X_{i_k,j_k},\ k=1,2$, for some $j_1\not=j_2$. 
The condition (4) in section 2 is shown to be equivalent to $(4)'$. 
Now also similar equivalence holds. 
Hence it suffices to see for each $x\in X'_{i'_2}\cap|T_{i'_1}^{X\prime}|$, the germ of $\pi_{i'_1}^{X\prime}|_{X'_{i'_2}\cap f^{-1}(f(x))}$ at $x$ is a $C^1$ submersion onto the germ of $X'_{i'_1}\cap f^{-1}(\pi_{j'_1}^{Y\prime}\circ f(x))$ at $\pi_{i'_1}^{X\prime}(x)$. 
We have four properties. 
$$
X'_{i'_2}\cap f^{-1}(f(x))=X_{i_2,j_2}\cap f^{-1}(f(x))\quad\text{by definition of }\{X'_{i'}\};
$$
$$
X'_{i'_1}\cap f^{-1}(\pi_{j'_1}^{Y\prime}\circ f(x))=X'_{i'_1}\cap f^{-1}(f\circ\pi_{i'_1}^{X\prime}(x))\quad\text{by (i)}
$$
$$
=X_{i_1,j_1}\cap f^{-1}(f\circ\pi_{i'_1}^{X\prime}(x))\quad\text{by definition of }\{X'_{i'}\};
$$
by $(4)'$ the germ of $\pi_{i_1,j_1}^X|_{X_{i_2,j_2}\cap f^{-1}(f(x))}$ at $x$ is a $C^1$ submersion onto the germ of $X_{i_1,j_1}\cap f^{-1}(f\circ\pi_{i_1,j_1}^X(x))$ at $\pi_{i_1,j_1}^X(x)$; by (iii) 
$$
\pi_{i'_1}^{X\prime}\circ\pi_{i_1,j_1}^X=\pi_{i'_1}^{X\prime}\quad\text{on}\ |T_{i'_1}^{X\prime}|\cap|T_{i_1,j_1}^X|. 
$$
Hence we only need to see the germ of $\pi_{i'_1}^{X\prime}|_{X_{i_1,j_1}\cap f^{-1}(f\circ\pi_{i_1,j_1}^X(x))}$ at $\pi_{i_1,j_1}^X(x)$ is a $C^1$ submersion onto the germ of $X'_{i'_1}\cap f^{-1}(f\circ\pi_{i'_1}^{X\prime}(x))$ at $\pi_{i'_1}^{X\prime}(x)$. 
That is clear by (i) because $f|_{X_{i_1,j_1}}:X_{i_1,j_1}\to Y_{j_1}$ is a $C^1$ submersion onto a union of some connected components of $Y_{j_1}$ and $f\circ\pi_{i_1,j_1}^X(x)$ and $f\circ\pi_{i'_1}^{X\prime}(x)$ are contained in the same connected component. \par
  Note we use the hypothesis $X'_{i'_k}\subset X_{i_k,j_k},\ k=1,2$, $j_1\not=j_2$ in the above arguments for only the property that the germ of $\pi_{i_1,j_1}^X|_{X_{i_2,j_2}\cap f^{-1}(f(x))}$ is a $C^1$ submersion into $X_{i_1,j_1}\cap f^{-1}(f\circ\pi_{i_1,j_1}^X(x))$. \par
  Case where $i_1\not=i_2$ and $X'_{i'_k}\subset X_{i_k,j_k},\ k=1,2$, for some $j_1$. 
In this case also the above property holds because $f\circ\pi_{i_1,j_1}^X=f$ on $X_{i_2,j_1}\cap|T_{i_1,j_1}^X|$ and $\pi_{i_1,j_1}^X|_{X_{i_2,j_1}\cap|T_{i_1,j_1}^X|}$ is a $C^1$ submersion into $X_{i_1,j_1}$. \par
  Case where $i_1=i_2$ and hence $X'_{i'_k}\subset X_{i_1,j_1},\ k=1,2$, for some $j_1$. 
In this case the reason is simply $\pi_{i_1.j_1}^X|_{X_{i_1,j_1}}=\id$. \par
  Thus (v) is proved. 
Now we define $\{T_{i'}^{X\prime}:i'\in I_{i_0}\}$. 
For that it suffices to consider separately $\{X'_{i'}:X'_{i'}\subset X_{i_0,j}\}$ for each $j$. 
Hence we assume all $X'_{i'}$ with $i'\in I_{i_0}$ are included in one $X_{i_0.j_0}$ for some $j_0$ and, moreover, $f(X_{i_0,j_0})=Y_{j_0}$ for simplicity of notation. 
Then as shown below we have a semialgebraic $C^1$ tube system $\{T_{i'}^{X\prime}=(|T_{i'}^{X\prime}|,\pi_{i'}^{X\prime},\rho_{i'}^{X\prime}):i'\in I_0\}$ for $\{X'_{i'}:i'\in I_0\}$ such that\par\noindent
(vi) for $i'$ and $j'$ with $i'\in I_{i_0}$ and $f(X'_{i'})=Y'_{j'}$, 
$$
f\circ\pi_{i'}^{X\prime}=\pi_{j'}^{Y\prime}\circ p\quad\text{on}\ |T_{i'}^{X\prime}|\cap p^{-1}(|T_{j'}^{Y\prime}|);
$$
(vii) for $j'\in J_{j_0},$ $\{T_{i'}^{X\prime}:f(X'_{i'})=Y'_{j'},\,i'\in I_{i_1},\,i_1\le i_0\}$ is a controlled semialgebraic $C^1$ tube system for $\{X'_{i'}:f(X'_{i'})=Y'_{j'},\,i'\in I_{i_1},\,i_1\le i_0\}$;\par\noindent
(viii) for $i_1,\,i'_k,\,k=1,2,3,\,i_4$ and $j_4$ with $i_1\le i_0,\,i'_1\in I_{i_1},\,i'_2,i'_3\in I_{i_0},\ X'_{i'_1}\cap(\overline{X'_{i'_2}}-X'_{i_2})\not=\emptyset$ and $X'_{i'_3}\subset\overline{X_{i_4,j_4}}$, 
$$
\pi_{i'_1}^{X\prime}\circ\pi_{i'_2}^{X\prime}=\pi_{i'_1}^{X\prime}\quad\text{on}\ |T_{i'_1}^{X\prime}|\cap|T_{i'_2}^{X\prime}|,
$$
$$
\pi_{i'_3}^{X\prime}\circ\pi_{i_4,j_4}^X=\pi_{i'_3}^{X\prime}\quad\text{on}\ |T_{i'_3}^{X\prime}|\cap|T_{i_4,j_4}^X|,
$$
if $i_0<i_4$ then 
$$
\rho_{i'_3}^{X\prime}\circ\pi_{i_4,j_4}^X=\rho_{i'_3}^{X\prime}\quad\text{on}\ |T_{i'_3}^{X\prime}|\cap|T_{i_4,j_4}^X|;
$$
(ix) for $i'\in I_{i_0}$ with $\dim X'_{i'}=i_0$, 
$$
T_{i'}^{X\prime}=T_{i_0,j_0}^X|_{|T_{i'}^{X\prime}|}. 
$$\par
  We construct $\{T_{i'}^{X\prime}:i'\in I_0\}$ as follows. 
First we define $T_{i'}^{X\prime}$ on $|T_{i'}^{X\prime}|\cap X_{i_0,j_0}$, $i'\in I_{i_0}$, so that (vi), (vii) and the first equality in (viii) are satisfied by the usual arguments of lift of a tube system (see [1], Lemma II.6.1, [4] and its proof). 
Secondly, extend $\pi_{i'}^{X\prime}$ to $|T_{i'}^{X\prime}|$ using $\pi_{i_0,j_0}^X$ as in the above construction of $\pi_{j'}^{Y\prime}$. 
Then $\pi_{i'}^{X\prime}$ are of class $C^1$; (vi) holds because for $i'$ and $j'$ with $i'\in I_{i_0}$ and $f(X'_{i'})=Y'_{j'}$, 
$$
f\circ\pi_{i'}^{X\prime}\overset{\text{definition of }\pi_{i'}^{X\prime}}=f\circ\pi_{i'}^{X\prime}\circ\pi_{i_0,j_0}^X\overset{\text{(vi) on }|T_{i'}^{X\prime}|\cap X_{i_0,j_0}}=\pi_{j'}^{Y\prime}\circ f\circ\pi_{i_0,j_0}^X
$$
$$
\overset{\text{(1) in section 2}}=\pi_{j'}^{Y\prime}\circ\pi_{j_0}^Y\circ p\overset{(**)_Y}=\pi_{j'}^{Y\prime}\circ p\quad\text{on}\ |T_{i'}^{X\prime}|\cap p^{-1}(|T_{j'}^{Y\prime}|);
$$
the first equality in (viii) for $i_1=i_0$ follows from definition of the extension; that for $i_1<i_0$ does from the second equality in (iii); the second in (viii) does from definition of the extension and the equality $\pi_{i_0,j_0}^X\circ\pi_{i_4,j_4}^X=\pi_{i_0,j_0}^X$; trivially $\pi_{i'}^{X\prime}=\pi_{i_0,j_0}^X$ for $i'\in I_{i_0}$ with $\dim X'_{i'}=i_0$. 
Thirdly, extend $\rho_{i'}^{X\prime}$ to $|T_{i'}^{X\prime}|$ in the same way as $\rho_{j'}^{Y\prime}$. 
Then $\{T_{i'}^{X\prime}:i'\in I_{i_0}\}$ is a semialgebraic $C^1$ tube system for $\{X'_{i'}:i'\in I_{i_0}\}$; (vii) holds because for $i'_0$ and $i'_1\in I_{i_1}$ with $i'_0\in I_{i_0},\ i_1<i_0$ and $f(X'_{i'_0})=f(X'_{i'_1})$, 
$$
\rho_{i'_1}^{X\prime}\circ\pi_{i'_0}^{X\prime}=\rho_{i'_1}^{X\prime}\circ\pi_{i'_0}^{X\prime}\circ\pi_{i_0,j_0}^X\quad\text{by definition of }\pi_{i'_0}^{X\prime}
$$
$$
=\rho_{i'_1}^{X\prime}\circ\pi_{i_0,j_0}^X\quad\text{by (vii) on }X_{i_0,j_0}
$$
$$
=\rho_{i'_1}^{X\prime}\quad\text{by the third equality in (iii)};
$$
the extensions are chosen so that the third equality in (viii) and (ix) are satisfied, which completes construction of a semialgebraic $C^1$ tube system $\{T_{i'}^{X\prime}:i'\in I_{i_0}\}$ and hence by induction that of $\{T_{i'}^{X\prime}:X'_{i'}\subset X\}$ with (i), (ii), the first equality in (iii) and (v) for any $i_0$, i.e., controlled over $\{T_{j'}^{Y\prime}\}$. \par
  It remains only to consider $X'_{i'}$ in $Z$, i.e., the case where $X'_{i'}$ is of the form $\{0\}\times Y'_{j'}$ for some $j'$. 
Set $\partial I=\{i':X'_{i'}\subset Z\}$. 
Obviously, we set 
$$
\pi_{i'}^{X\prime}(x)=(0,\pi_{j'}^{Y\prime}\circ p(x))\quad\text{for}\ x\in|T_{i'}^{X\prime}|,\ i'\in\partial I\ \text{and}\ j'\ \text{with}\ X'_{i'}=\{0\}\times Y'_{j'},
$$
where $|T_{i'}^{X\prime}|$ is a small semialgebraic neighborhood of $X'_{i'}$ in $\R^m\times\R^n$. 
Then (i) for $i'\in\partial I$ is clear; the first equality in (iii) for $i'_1\in\partial I$ holds because 
$$
\pi_{i'_1}^{X\prime}\circ\pi_{i'_2}^{X\prime}(x)\overset{\text{definition of }\pi_{i'_1}^{X\prime}}=(0,\pi_{j'_1}^{Y\prime}\circ p\circ\pi_{i'_2}^{X\prime}(x))\overset{\text{(i)}}=(0,\pi_{j'_1}^{Y\prime}\circ\pi_{j'_2}^{Y\prime}\circ p(x))
$$
$$
\overset{\text{controlledness of }\{T_{j'}^{Y\prime}\}}=(0,\pi_{j'_1}^{Y\prime}\circ p(x))=\pi_{i'_1}^{X\prime}(x)\quad\text{for}\ x\in|T_{i'_1}^{X\prime}|\cap|T_{i'_2}^{X\prime}|, 
$$
where $j'_1$ and $j'_2$ are such that $f(X'_{i'_k})=Y'_{j'_k},\ k=1,2$; (v) for $i'_1\in\partial I$ is clear, to be precise, for $i'_1\in\partial I,\ i'_2,\ j'_1$ and $j'_2$ with $X'_{i'_1}\cap(\overline{X'_{i_2}}-X'_{i'_2})\not=\emptyset$, $Y'_{j'_1}\subset\overline{Y'_{j'_2}}-Y'_{j'_2}$ and $p(X'_{i'_k})=Y'_{j'_k},\ k=1,2$, the restriction of the map $(\pi_{i'_1}^{X\prime},p)$ to $X'_{i'_2}\cap|T_{i'_1}^{X\prime}|$ is a $C^1$ submersion into $X'_{i'_1}\times_{(p,\pi_{j'_1}^{Y\prime})}(Y'_{j'_2}\cap|T_{j'_1}^{Y\prime}|)$ because $p|_{X'_{i'_1}}:X'_{i'_1}\to Y'_{j'_1}$ is a $C^1$ diffeomorphism and $p|_{X'_{i'_2}}:X'_{i'_2}\to Y'_{j'_2}$ is a $C^1$ submersion. \par
  We want to define $\{\rho_{i'}^{X\prime}:i'\in\partial I\}$ so that $\{T_{i'}^{X\prime}=(|T_{i'}^{X\prime}|,\pi_{i'}^{X\prime},\rho_{i'}^{X\prime}):i'\in\partial I\}$ is a semialgebraic $C^1$ weak tube system and for each $j'$, $\{T_{i'}^{X\prime}:f(X'_{i'})=Y'_{j'}\}$ is controlled. 
We proceed by double induction. 
Let $d\ge 0\in\Z$, and assume $\rho_{i'}^{X\prime}$ are already defined if $\dim X'_{i'}>d$. 
We need to construct $\rho_{i'}^{X\prime}$ for $i'\in\partial I$ with $\dim X'_{i'}=d$. 
As the problem is local at such $X'_{i'}$, assume there exists only one $i'_0\in\partial I$ with $\dim X'_{i'_0}=d$. 
Set $I'=\{i':X'_{i'_0}\subset\overline{X'_{i'}}-X'_{i'}\}$ and $Y'_{j'_0}=p(X'_{i'_0})$. \par 
  For the moment we construct a non-negative semialgebraic $C^0$ function $\rho_{i'_0,d}^{X\prime}$ on $|T_{i'_0}^{X\prime}|$ with zero set $X'_{i'_0}$ which is of class $C^1$ on $|T_{i'_0}^{X\prime}|-X'_{i'_0}$ and such that $\{T_{i'_0,d}^{X\prime},T_{i'}^{X\prime}:i'\in I',\,p(X'_{i'})=Y'_{j'_0}\}$ is controlled, i.e., 
$$
\rho_{i'_0,d}^{X\prime}\circ\pi_{i'}^{X\prime}=\rho_{i'_0,d}^{X\prime}\quad\text{on}\ |T_{i'_0}^{X\prime}|\cap|T_{i'}^{X\prime}|\ \text{for}\ i'\in I'\ \text{with}\ p(X'_{i'})=Y'_{j'_0}, 
$$
where $d=1+\#I'$ and $T_{i'_0,d}^{X\prime}=(|T_{i'_0}^{X\prime}|,\pi_{i'_0}^{X\prime},\rho_{i'_0,d}^{X\prime})$. 
(Namely we forget the condition that $\rho_{i'_0,d}^{X\prime}|_{X'_{i'}\cap\pi_{i'_0}^{X\prime-1}(x)-X'_{i'_0}}$ is $C^1$ regular for each $x$ and any $i'\in I'$.) 
Order elements of $I'$ as $\{i'_1,...,i'_{d-1}\}$ so that $\dim X'_{i'_1}\le\cdots\le\dim X'_{i'_{d-1}}$. \par
  Let $k\in\Z$ with $0\le k<d-1$. 
As the second induction, assume we have a non-negative semialgebraic $C^0$ function $\rho_{i'_0,k}^{X\prime}$ defined on $|T_{i'_0}^{X\prime}|\cap(|T_{i'_1}^{X\prime}|\cup\cdots\cup|T_{i'_k}^{X\prime}|)$ such that $\rho_{i'_0,k}^{X\prime-1}(0)=X'_{i'_0},$ $\rho_{i'_0,k}^{X\prime}$ is of class $C^1$ outside of $X'_{i'_0}$ and $\{T_{i'_0,k}^{X\prime},T_{i'}^{X\prime}:i'\in I',\,p(X'_{i'})=Y'_{j'_0}\}$ is controlled, i.e., 
$$
\rho_{i'_0,k}^{X\prime}\circ\pi_{i'}^{X\prime}=\rho_{i'_0,k}^{X\prime}\quad\text{on}\ |T_{i'_0}^{X\prime}|\cap(|T_{i'_1}^{X\prime}|\cup\cdots\cup|T_{i'_k}^{X\prime}|)\cap|T_{i'}^{X\prime}|\ \text{for}\ i'\in I'\ \text{with}\ p(X'_{i'})=Y'_{j'_0}, 
$$
where $T_{i'_0,k}^{X\prime}=(|T_{i'_0}^{X\prime}|,\pi_{i'_0}^{X\prime},\rho_{i'_0,k}^{X\prime})$. 
Then we need to define $\rho_{i'_0,k+1}^{X\prime}$. 
Let $\tilde\rho_{i'_0,k}^{X\prime}$ be any non-negative semialgebraic $C^0$ extension of $\rho_{i'_0,k}^{X\prime}|_{|T_{i'_0}^{X\prime}|\cap(|T_{i'_1}^{X\prime}|\cup\cdots\cup|T_{i'_k}^{X\prime}|)\cap X'_{i'_{k+1}}}$ to $|T_{i'_0}^{X\prime}|\cap X'_{i'_{k+1}}$ with zero set $X'_{i'_0}$, let $V$ be an open semialgebraic neighborhood of $X'_{i'_1}\cup\cdots\cup X'_{i'_k}$ in $X'_{i'_1}\cup\cdots\cup X'_{i'_{k+1}}$ whose closure is included in $|T_{i'_1}^{X\prime}|\cup\cdots\cup |T_{i'_k}^{X\prime}|$, approximate $\tilde\rho_{i'_0,k}^{X\prime}|_{|T_{i'_0}^{X\prime}|\cap X'_{i'_{k+1}}-V}$ by a non-negative semialgebraic $C^0$ function $\Tilde{\Tilde\rho}_{i'_0,k}^{X\prime}$ in the uniform $C^0$ topology so that $\Tilde{\Tilde\rho}_{i'_0,k}^{X\prime-1}(0)=X'_{i'_0}$, and $\Tilde{\Tilde\rho}_{i'_0,k}^{X\prime}$ is of class $C^1$ outside of $X'_{i'_0}$ (Theorem II.4.1, [3]), let $\xi$ be a semialgebraic $C^1$ function on $|T_{i'_0}^{X\prime}|\cap X'_{i'_{k+1}}$ such that $0\le\xi\le 1$, $\xi=0$ on $|T_{i'_0}^{X\prime}|\cap X'_{i'_{k+1}}\cap V$ and $\xi=1$ on $|T_{i'_0}^{X\prime}|\cap X'_{i'_{k+1}}-|T_{i'_1}^{X\prime}|-\cdots|T_{i'_k}^{X\prime}|$, and set 
$$
\hat\rho_{i'_0,k}^{X\prime}(x)=\xi(x)\Tilde{\Tilde\rho}_{i'_0,k}^{X\prime}(x)+(1-\xi(x))\rho_{i'_0,k}^{X\prime}(x)\quad\text{for }x\in|T_{i'_0}^{X\prime}|\cap X'_{i'_{k+1}}. 
$$
Then $\hat\rho_{i'_0,k}^{X\prime}$ is a non-negative semialgebraic $C^0$ extension of $\rho_{i'_0,k}^{X\prime}|_{|T_{i'_0}^{X\prime}|\cap V\cap X'_{i'_{k+1}}}$ to $|T_{i'_0}^{X\prime}|\cap X'_{i'_{k+1}}$ with zero set $X'_{i'_0}$ and of class $C^1$ outside of $X'_{i'_0}$. 
If $p(X'_{i'_{k+1}})\not=Y'_{j'_0}$, we continue to extend $\hat\rho_{i'_0,k}^{X\prime}$ to the required $\hat\rho_{i'_0,k+1}^{X\prime}:|T_{i'_0}^{X\prime}|\cap(|T_{i'_1}^{X\prime}|\cup\cdots\cup|T_{i'_{k+1}}^{X\prime}|)\to\R$ shrinking $|T_{i'_1}^{X\prime}|,...,|T_{i'_k}^{X\prime}|$ and using a partition of unity in the same way so that $\rho_{i'_0,k+1}^{X\prime}=\rho_{i'_0,k}^{X\prime}$ on $|T_{i'_0}^{X\prime}|\cap(|T_{i'_1}^{X\prime}|\cup\cdots\cup|T_{i'_k}^{X\prime}|)$. 
Otherwise, set 
$$
\rho_{i'_0,k+1}^{X\prime}=\left\{
\begin{array}{ll}
\rho_{i'_0,k}^{X\prime}\quad&\text{on}\ |T_{i'_0}^{X\prime}|\cap(|T_{i'_1}^{X\prime}|\cup\cdots\cup|T_{i'_k}^{X\prime}|)\\
\hat\rho_{i'_0,k}^{X\prime}\circ\pi_{i'_0,k+1}^{X\prime}\quad&\text{on}\ |T_{i'_0}^{X\prime}|\cap|T_{i'_{k+1}}^{X\prime}|, \\
\end{array}
\right.
$$
which is well-defined because 
$$
\hat\rho_{i'_0,k}^{X\prime}\circ\pi_{i'_0,k+1}^{X\prime}=\rho_{i'_0,k}^{X\prime}\circ\pi_{i'_0,k+1}^{X\prime}\quad\text{by definition of }\hat\rho_{i'_0,k}^{X\prime}\qquad\qquad\qquad\quad
$$
$$
=\rho_{i'_0,k}^{X\prime}\quad\text{by controlledness of }\{T_{i'_0,k}^{X\prime},T_{i'}^{X\prime}:i'\in I',\ p(X'_{i'})=Y'_{j'_0}\}\qquad
$$
$$
\qquad\qquad\text{on}\ |T_{i'_0}^{X\prime}|\cap(|T_{i'_1}^{X\prime}|\cup\cdots\cup|T_{i'_k}^{X\prime}|)\cap|T_{i'_{k+1}}^{X\prime}|\text{ for shrunk }|T_{i'_1}^{X\prime}|,...,|T_{i'_k}^{X\prime}|. 
$$
Then clearly $\rho_{i'_0,k+1}^{X\prime-1}(0)=X'_{i'_0},\ \rho_{i'_0,k+1}^{X\prime}$ is of class $C^1$ outside of $X'_{i'_0}$ and 
$$
\rho_{i'_0,k+1}^{X\prime}\circ\pi_{i'}^{X\prime}=\rho_{i'_0,k+1}^{X\prime}\ \text{on}\ |T_{i'_0}^{X\prime}|\cap(|T_{i'_1}^{X\prime}|\cup\cdots\cup|T_{i'_{k+1}}^{X\prime}|)\cap|T_{i'}^{X\prime}|\ \text{for}\ i'\in I'\ \text{with}\ p(X'_{i'})=Y'_{j'_0}
$$
as follows. It suffices to consider only the case where $\overline {X'_{i'}}-X'_{i'}\supset X'_{i'_{k+1}}$ and $p(X'_{i'})=p(X'_{i'_{k+1}})=Y'_{j'_0}$ and the equation on $|T_{i'_0}^{X\prime}|\cap|T_{i'_{k+1}}^{X\prime}|\cap|T_{i'}^{X\prime}|$. 
We have 
$$
\rho_{i'_0,k+1}^{X\prime}\circ\pi_{i'}^{X\prime}=\hat\rho_{i'_0,k}^{X\prime}\circ\pi_{i'_{k+1}}^{X\prime}\circ\pi_{i'}^{X\prime}\quad\text{by definition of }\rho_{i'_0,k+1}^{X\prime}\qquad\qquad
$$
$$
=\hat\rho_{i'_0,k}^{X\prime}\circ\pi_{i'_{k+1}}^{X\prime}\quad\text{by the first equation in (iii)}
$$
$$
\qquad\qquad=\rho_{i'_0,k+1}^{X\prime}\quad\text{by definition of }\rho_{i'_0,k+1}^{X\prime}\quad\text{on}\ |T_{i'_0}^{X\prime}|\cap|T_{i'_{k+1}}^{X\prime}|\cap|T_{i'}^{X\prime}|. 
$$
Thus by the second induction we obtain $\rho_{i'_0,d-1}^{X\prime}:|T_{i'_0}^{X\prime}|\cap(|T_{i'_1}^{X\prime}|\cup\cdots\cup|T_{i',d-1}^{X\prime}|)\to\R$. \par
  It remains only to extend $\rho_{i'_0,d-1}^{X\prime}$ to a non-negative semialgebraic $C^0$ function $\rho_{i'_0,d}^{X\prime}$ on $|T_{i'_0}^{X\prime}|$ with zero set $X'_{i'_0}$ and of class $C^1$ outside of $X'_{i'_0}$. 
However we have already carried out such a sort of extension by using a partition of unity $\xi$. \par
  We need to solve the problem of $C^1$ regularity of $\rho_{i'_0,d}^{X\prime}|_{X'_{i'}\cap\pi_{i'_0}^{X\prime-1}(x)-X'_{i'_0}}$. 
For each $x\in X'_{i'_0}$, the restriction of $\rho_{i'_0,d}^{X\prime}$ to $X'_{i'}\cap\pi_{i'_0}^{X\prime-1}(x)\cap\rho_{i'_0,d}^{X\prime-1}((0,\,\delta_x))$ is $C^1$ regular for some $\delta_x>0\in\R$ and any $i'\in I'$. 
Here we can choose $\delta_x$ so that the function $X'_{i'_0}\ni x\to\delta_x\in\R$ is semialgebraic (but not necessarily continuous). 
Then there exists a semialgebraic closed subset $X''_{i'_0}$ of $X'_{i'_0}$ of smaller dimension such that each point $x$ in $X'_{i'_0}-X''_{i'_0}$ has a neighborhood in $X'_{i'_0}$ where $\delta_x$ is larger than a positive number. 
Hence if we replace $X'_{i'_0}$ with $X'_{i'_0}-X''_{i'_0}$, i.e., $Y'_{j'_0}$ with $Y'_{j'_0}-p(X''_{i'_0})$ and shrink $|T_{i'_0}^{X\prime}|$ then the $C^1$ regularity holds. 
Thus we obtain the required $\rho_{i'_0}^{X\prime}$ though $X'_{i'_0}$ is shrunk to $X'_{i'_0}-X''_{i'_0}$. \par
  The shrinking is admissible as follows. 
Substratify $\{Y'_{j'}\cap p(\overline{X''_{i'_0}}),\,Y'_{j'}-p(\overline{X''_{i'_0}})\}$ to a Whitney semialgebraic $C^1$ stratification $\{Y''_{j''}\}$ such that $\{Y'_{j'}-p(\overline{X''_{i'_0}})\}=\{Y''_{j''}-p(\overline{X''_{i'_0}})\}$, set $\{X''_{i''}\}=\{X_{i,j}\cap p^{-1}(Y''_{j''}),\,Z\cap\{0\}\times Y''_{j''}\}$, which implies $\{X'_{i'}-p^{-1}(p(\overline{X''_{i'_0}}))\}=\{X''_{i''}-p^{-1}(p(\overline{X''_{i'_0}}))\}$, and repeat all the above arguments to $\overline f:\{X''_{i''}\}\to\{Y''_{j''}\}$. 
Then we obtain a semialgebraic $C^1$ tube system $\{T_{j''}^{Y\prime\prime}\}$ for $\{Y''_{j''}\}$ and a semialgebraic $C^1$ tube system $\{T_{i''}^{X\prime\prime}:X''_{i''}\subset X\}$ for $\{X''_{i''}\subset X\}$ controlled over $\{T_{j''}^{Y\prime\prime}\}$ such that $\{T_{j''}^{Y\prime\prime}:Y''_{j''}\cap p(\overline{X''_{i'_0}})=\emptyset\}$ and $\{T_{i''}^{X\prime\prime}:X''_{i''}\subset X,\ X''_{i''}\cap p^{-1}(p(\overline{X''_{i'_0}}))=\emptyset\}$ are equal to $\{T_{j'}^{Y\prime}|_{|T_{j'}^{Y\prime}|-\pi_{j'}^{Y\prime-1}(p(\overline{X''_{i'_0}}))}\}$ and $\{T_{i'}^{X\prime}|_{|T_{i'}^{X\prime}|-\pi_{i'}^{X\prime-1}(p^{-1}(p(\overline{X''_{i'_0}})))}\}$, respectively, by (iv) and (ix), where the domains of the latter two tube systems are shrunk. 
Moreover we continue construction of $\rho_{i''}^{X\prime\prime}$ for $X''\subset Z$. 
Since $\{X''_{i''}\subset Z:\dim X''_{i''}>d\}=\{X'_{i'}\subset Z:\dim X'_{i'}>d\}$ and $\{X''_{i''}\subset Z:\dim X''_{i''}=d\}=\{X'_{i'_0}-X''_{i'_0}\}$ we choose $\rho_{i'}^{X\prime}$ as $\rho_{i''}^{X\prime\prime}$ for $X''_{i''}\subset Z$ with $\dim X''_{i''}>d$ and $\rho_{i'_0}^{X\prime}|_{|T_{i''}^{X\prime\prime}|}$ as $\rho_{i''}^{X\prime\prime}$ for $X''_{i''}\subset Z$ with $\dim X''_{i''}=d$. 
Hence we can assume $X''_{i'_0}=\emptyset$ from the beginning, which completes the construction of $\rho_{i'_0}^{X\prime}$ and hence of the required $\{\rho_{i'}^{X\prime}:i'\in\partial I\}$ by induction.\par
  Thus $\overline f:\{X'_{i'}\}\to\{Y'_{j'}\},\ \{T_{i'}^{X\prime}\}$ and $\{T_{j'}^{Y\prime}\}$ satisfy the conditions in theorem 2.2. 
Hence theorem 1.2 follows. \qed

\end{document}